\documentclass[10pt]{amsart}\usepackage{amssymb} \usepackage{amscd}

\usepackage{amssymb}
\usepackage{amscd}
\usepackage{epsfig}
\usepackage[all]{xy}

\newcommand{\Q}{\mathbb Q}
\newcommand{\Z}{\mathbb Z}
\newcommand{\C}{\mathbb C}
\newcommand{\cc}{\mathcal C}

\newcommand{\h}{\mbox{H}}
\newcommand{\N}{\mathcal N}
\newcommand{\hm}{\mbox{Hom}}
\newcommand{\I}{\mathbb I}
\newcommand{\D}{\mathcal D}
\newcommand{\To}{\longrightarrow}
\newtheorem{thm}{Theorem}[section]
\newtheorem{lem}[thm]{Lemma}

\newtheorem{prop}[thm]{Proposition}

{\theoremstyle{definition} \newtheorem{rem}[thm]{Remark}
} {\theoremstyle{definition}
\newtheorem{defn}[thm]{Definition}}

{\theoremstyle{remark} }
\title[Motivic Iterated Integrals]{Algebraic Cycles and Motivic Iterated
Integrals II}
\author{Amir Jafari}

\begin{document}

\bibliographystyle{amsalpha+}

\begin{abstract}
 This is a sequel to \cite{FJ}. It will give a more natural framework for constructing elements in the Hopf algebra
 $\chi_F$ of framed mixed Tate motives according to Bloch and K\v ri\v z \cite{BK}. This framework allows us to extend the
 results of \cite{FJ} to interpret all multiple zeta values (including the divergent ones) and the multiple polylogarithms in
 one variable as elements of $\chi_F$. It implies that the pro-unipotent completion of the torsor of paths on ${\mathbb P}^1-\{0,1,\infty\}$,
 is a mixed Tate motive in the sense of \cite{BK}. Also It allows us to interpret the multiple logarithm $Li_{1,\dots,1}(z_1,\dots,z_n)$ as an
 element of $\chi_F$ as long as the products of consecutive $z_i$'s are not $1$.
\end{abstract}

\maketitle
\tableofcontents


\section{Introduction}
The category $\mathcal{MTM}(F)$ of mixed Tate motives over a field $F$ is supposed to be a Tannakian category with a special invertible object $\Q(1)$. Let $\Q(-1)=\underline{\hm}(\Q,\Q(1))$. The simple objects of this category are the mutually distinct elements $\Q(n)=\otimes^{|n|}\Q(n/|n|)$ and $\Q(0)=\Q$. If we define the weight of $\Q(n)$ to be $-2n$, then any object of this category has a canonical increasing (weight) filtration such that $gr_{-2n+1}=0$ and $gr_{-2n}$ is a direct sum of $\Q(n)$'s. One should also have the following relation with $K$-theory of $F$:
$$\mbox{Ext}^i(\Q(0),\Q(n))=gr_{\gamma}^n\mbox{K}_{2n-i}(F)_{\Q}.$$
Furthermore if $F$ can be embedded inside the field of complex numbers $\C$, there should be a Hodge realization tensor functor from $\mathcal{MTM}(F)$ to the category of mixed Hodge-Tate structures. $\Q(1)$ will correspond to the Hodge structure of $\h_2({\mathbb P}^1)$ or equivalently $\h_1({\mathbb P}^1\backslash\{0,\infty\})$. Also when $F$ has a characteristic prime to a prime number $\ell$ and contains only finitely many $\ell^{th}$ power roots of $1$, there should be an $\ell$-adic realization tensor functor from $\mathcal{MTM}(F)$ to the category of mixed Tate $\ell$-adic representations of the absolute Galois group $G_F = \mbox{Gal}(\overline{F}/F)$. $\Q(1)$ will correspond to the action of $G_F$ on
$$\Q \otimes \lim_{\leftarrow}\mu_{\ell^n}(\bar{F})=\Q_{\ell}(1).$$

The existence of this category in this generality is still conjectural. However if the field $F$ satisfies Beilinson-Soul\' e conjecture stating:
$$gr_{\gamma}^n\mbox{K}_{2n-i}(F)_{\Q}=0\quad \mbox{for}\quad i\le 0 \quad \mbox{(except when $i=n=0$)}$$
such a category can be constructed as the heart of a canonical $t$-structure on the triangulated subcategory of mixed motives over $F$ (constructed independently by Voevodsky, Levine and Hanamura) generated by $\Q(n)$'s. It is known for example thanks to Borel that a number field satisfies this conjecture.

Bloch and K\v ri\v z in \cite{BK} proposed a different approach. Using algebraic cycles in $\Box_F^n = ({\Bbb P}^1_F\backslash \{1\})^n$ whose support meets all the faces (defined by setting some of the coordinates $= 0$ or $\infty$) properly, they constructed a graded Hopf algebra that they proposed to be the Hopf algebra of framed mixed Tate motives.

Let $\mbox{Cycle}^r(n)$ be admissible cycles of codimension $r$ in $\Box_F^n$.
The permutation group $\Sigma_n$ and the group $(\Z/2\Z)^n$ act on $\Box_F^n$. So we have an action of ${\mathcal G}_n = \Sigma_n\rtimes (\Z/2\Z)^n$ on $\mbox{Cycle}^r(n)$. Let
$$\N^n(r)=\mbox{Alt}(\mbox{Cycle}^r(2r-n)\otimes \Q)$$
where
$$\mbox{Alt}=\frac{1}{|\mathcal G_{2r-n}|}\sum_{\tau\in {\mathcal G}_{2r-n}} \mbox{sgn}(\tau)\tau.$$
The sign of an element $(\sigma, \epsilon_1,\dots,\epsilon_n)$ where $\sigma\in \Sigma_n$ and $\epsilon_i=\pm 1$ is defined to be $\mbox{sgn}(\sigma)\epsilon_1\dots\epsilon_n$. The external product followed by Alt defines a product
$$\N^n(r)\otimes \N^m(s)\To \N^{n+m}(r+s)$$
that is graded commutative with respect to the upper grading. Also we have a differential:
$$d:\N^n(r)\To \N^{n+1}(r)$$
$$Z\mapsto \sum_{i=1}^{2r-n} (-1)^{i-1}(Z|_{x_i=0}-Z|_{x_i=\infty}).$$
These structures make $\N=\oplus\N^n(r)$ into a graded differential algebra (DGA) with Adams grading. The seminal works of Bloch and Levine on higher Chow groups imply that:
$$\h^i(\N(n)^{\bullet})=gr_{\gamma}^n \mbox{K}_{2n-i}(F)_{\Q}.$$
Hence the conjecture of Beilinson-Soul\'e briefed above is equivalent to the fact that $\N$ is cohomologically connected.

The bar construction, invented by Eilenberg-MacLane and Adams, to any DGA  $\mathcal A$ associates a new DGA $B(\mathcal A)$ with a coproduct. This construction is reminded in appendix A.

The proposed Hopf algebra of framed mixed Tate motives is then
$$\chi_F = \h^0B(\N_F).$$

In \cite{BK} the realization morphisms from $\chi_F$ to the Hopf algebras of framed mixed Hoge-Tate structures and framed mixed $l$-adic representations of $G_F$ are constructed.

The category of mixed Tate motives over $F$ is defined to be the graded co-modules over $\chi_F$.

A rich source for constructing mixed motives is the unipotent completion of the fundamental group or the torsor of paths of smooth algebraic varieties. For any smooth variety $X$ and $a,b\in X$, the truncated torsor of paths  $\Pi^{N}(X,a,b)$ is defined by:
$$\frac{\Q[\pi(X,a,b)]}{I_a^{N+1}\Q[\pi(X,a,b)]}.$$
Here $I_a$ is the augmentation ideal $\mbox{Ker}(\Q[\pi_1(X,a)]\To \Q)$ (We are viewing $X$ as a complex manifold by taking its points over $\C$, however there is an algebraic way to define these objects). According to Morgan and Hain this space has a natural mixed Hodge structure, in fact Hain and Zucker showed that $\Pi^{N}(X,a,b)$ defines a good unipotent variation of mixed Hodge structures on $X\times X$. Therefore even for $a$ and $b$ in a smooth compactification $\bar{X}$ of $X$ and not in $X$, one can take limit mixed Hodge structure to define the truncated torsor of paths from $a$ to $b$. It's easy to see that $gr_0\Pi^N=\Q(0)$. Assuming $W_{-4}H_2(X)=0$, one can show
$$gr_{-2N}\Pi^N(X,a,b)=\otimes^N W_{-2}H_1(X)$$
(Here equality means a canonical isomorphism.)
Therefore if one has non-vanishing holomorphic functions $f_1,\dots,f_N$ on $X$, we get a morphism:
$$f_{1*}\otimes \dots\otimes f_{N*}: \otimes^N W_{-2}H_1(X)\To \otimes^N W_{-2}H_1({\Bbb P}^1\backslash\{0,\infty\})=\Q(N).$$

Hence we get a framed mixed Hodge structure with the canonical frame $gr_0\Pi^N =\Q(0)$ and the frame:
$$f_{1*}\otimes \dots\otimes f_{N*}:gr_{-2N}\Pi^N\To \Q(N).$$
This is the Hodge theoretic version of the Iterated integral:
$$\int_a^b \frac{df_1}{f_1}\circ \dots \circ \frac{df_N}{f_N}.$$

For convenience of the reader we recall here the definition of such integrals. Let $M$ be a manifold and $\omega_1,\dots, \omega_n$ smooth 1-forms on $M$. For any path $\gamma:[0,1]\To M$, let $\gamma^*\omega_i =f_i(t)dt$, then:
$$\int_{\gamma}\omega_1\circ \dots\circ \omega_n =\int_{0\le t_1\le\dots\le t_n\le 1} f_1(t_1)\dots f_n(t_n)dt_1\dots dt_n.$$

The framed mixed Hodge structure defined above, being of geometric origin, should  have an upgrade to a framed mixed motive. If one assumes that for a smooth compactification $\bar{X}$ of $X$, $\h^1(\bar{X})=0$, then it is shown that $\Pi^N(X,a,b)$ carries a mixed Hodge-Tate structure and therefore we expect to have an upgrade to a framed mixed Tate motive and therefore an element ${\Bbb I}_X(a;f_1,\dots,f_N;b) \in \h^0B(\N)$.

Our goal in this article is to carry out this problem for the case $X={\Bbb P}^1\backslash\{0,1,\infty\}$ completely and therefore giving motivic interpretation of the multiple zeta values:
$$\zeta(n_1,\dots,n_k)=\sum_{0<m_1<\dots<m_k}\frac{1}{m_1^{n_1}\dots m_k^{n_k}}$$
and more generally a motivic interpretation of the multiple polylogarithms in one variable:
$$Li_{n_1,\dots,n_k}(z)=\sum_{0<m_1<\dots<m_k}\frac{z^{m_k}}{m_1^{n_1}\dots m_k^{n_k}}.$$
As they can be written as iterated integrals on ${\mathbb P}^1\backslash\{0,1,\infty\}$:
$$\int_0^z \underbrace{\frac{dt}{1-t}\circ \frac{dt}{t}\circ\dots\circ \frac{dt}{t}}_{n_1}\circ \dots\circ \underbrace{\frac{dt}{1-t}\circ\frac{dt}{t}\circ\dots \circ \frac{dt}{t}}_{n_k}.$$

Unfortunately since some of the cycle constructions are not admissible we are unable to treat the more general case
for the multi-variable multiple polylogarithms:
$$Li_{n_1,\dots,n_k}(z_1,\dots,z_k)=\sum_{0<m_1<\dots<m_k}\frac{z_1^{m_1}\dots z_k^{m_k}}{m_1^{n_1}\dots m_k^{n_k}}$$
which also has an iterated integral definition (assuming $z_i\ne 0$):
$$\int_0^1 \underbrace{\frac{dt}{a_1-t}\circ \frac{dt}{t}\circ\dots\circ \frac{dt}{t}}_{n_1}\circ \dots\circ \underbrace{\frac{dt}{a_k-t}\circ\frac{dt}{t}\circ\dots \circ \frac{dt}{t}}_{n_k}.$$
with $a_i=(z_iz_{i+1}\dots z_k)^{-1}$.
 However we can give a motivic interpretation of the iterated integral:
$$\int_{a_0}^{a_{n+1}}\frac{dt}{t-a_1}\circ\dots\circ\frac{dt}{t-a_n}$$
in the following cases:
\begin{enumerate}
\item There is an element $a$ in $F$ such that $a_i=a_j$ for $i\ne j$ implies $a_i=a$. We call such a sequence $a$-generic.

\item For $i=1,\dots,n$, we have $a_i\ne a_{i+1}$. We call such a sequence sequentially distinct.
\item There are two distinct elements $\alpha$ and $\beta$ such that $a_i\in\{\alpha,\beta\}$ for $i=1,\dots,n$. We call such a sequence binary.
\end{enumerate}

For example the second condition allows us to give an element in $\h^0B(\N)$ that corresponds to the multiple logarithm $Li_{1,\dots,1}(z_1,\dots,z_k)$ with the assumption that for any $m\le n$: $\prod_{i=m}^n z_i\ne 1$. With the more restrictive assumption on $z_i$'s which demands $\prod_{i\in I} z_i\ne 1$ for any subset $I$ of $\{1,\dots,k\}$ the first condition can be used to construct an element in $\h^0B(\N)$ that corresponds to the multi-variable multiple polylogarithm $Li_{n_1,\dots,n_k}(z_1,\dots,z_k)$.
 \begin{rem}
 In our joint work with Furusho, we treated the case of $0$-generic sequences in the above terminology. This article is written so that it will redo the results of \cite{FJ} in a more transparent framework and extend them to cover more cases. Therefore there are some major overlaps; the author has decided to keep them for the sake of convenience for  the reader. The language of trees, which was also used by Gangl, Goncharov and Levin in \cite{GGL}, has been replaced by the combinatorial DGA introduced in section 2. It has to be stressed that the existence of this paper and \cite{FJ} owe a great deal to \cite{GGL}.
\end{rem}
\section{A Combinatorial Differential Graded Algebra}

Let $S$ be a set and $\cc_S$  the graded commutative algebra with unit, freely generated by symbols:
$$(a_0;a_1,a_2,\dots,a_n;a_{n+1})$$
of degree 1 and Adams grading $n$ where $a_i\in S$ and $n>0$.
\begin{defn}
A 2-cut of $A=(a_0;a_1,\dots,a_n;a_{n+1})$ is an ordered pair
$$(a_0;\dots,a_{i-1},a_i,a_{j+1},a_{j+2},\dots;a_{n+1}), (a_i;a_{i+1},\dots;a_j;a_{j+1})$$ for a choice $0\le i< j \le n$ , $(i,j)\ne (0,n)$.

Similarly a $k$-cut is an ordered $k$-tuple obtained from a $(k-1)$-cut $A_1,\dots,A_{k-1}$ by replacing one of the $A_i$'s with a 2-cut $A_i'$ and $A_i''$. We also define a 1-cut of $A$ to be just $A$. Note that a cut should be regarded in terms of the indices and not the elements. So two 2-cuts that correspond to two different choices for $i<j$ are considered to be different even if their realized pairs are identical.
\end{defn}

Define the differential on $\cc_S$  for a generator $A=(a_0;a_1,\dots,a_n;a_{n+1})$ to be
$$dA = -\sum A'A''$$
where the sum is over all 2-cuts of $A$.
\begin{lem}
$d^2=0.$
\end{lem}
{\bf Proof.} Let $A' = (a_0,\dots,a_i,a_{j+1},\dots,a_{n+1}), A'' = (a_i;\dots;a_{j+1})$ be a 2-cut of $A$. A 2-cut $A_1',A_2'$ of $A'$ is said to be joint with $A''$ if $A_2',A''$ is a 2-cut of a subsequence $A_3=(a_k;a_{k+1},\dots;a_{l+1})$ of $A$. otherwise we call $(A_1',A_2'')$ disjoint from $A''$.
Now
$$ddA=\sum A_1'A_2'A''-\sum A'A_1''A_2''.$$
The terms in this sum cancel each other as follows. If $A_1',A_2'$ is disjoint from $A''$ then can consider removing $A_2'$ from $A$ first and then remove $A''$ from the remaining of $A$, we will get $A_1'A''A_2'$ that cancels $A_1'A_2'A''$ due to graded commutativity of $\cc_S$. If the cut is joint and $A_2',A''$ is a 2-cut of $A_3$ as above, we can remove $A_3$ from $A$ first and we get the term $A_1'A_3$  and then cut $A_3$ into $A_2'$ and $A''$, this cancels $A_1'A_2'A''$ since it appears with a negative sum in the formula above for $d^2$. A similar argument for terms of the form $A'A_1''A_2''$ finishes the proof.\qed
\begin{defn}
Associated to a generator $A=(a_0;a_1,\dots,a_n;a_{n+1})$ define $T(A)\in B(\cc_S)^0$ by:
$$T(A)=\sum [A_1|\dots|A_k]$$
where the sum is over all cuts of $A$ (including 1-cut).
\end{defn}
\begin{lem}
For $A=(a_0;a_1,\dots,a_n;a_{n+1})$, $d(T(A))=0$. Therefore $T(A)\in H^0B(\cc_S)(n)\subseteq B(\cc_S)^0$.
\end{lem}
{\bf Proof.} This follows once we observe:
$$d_{ext}(\sum_{\mbox{\small{k-cuts}}}[A_1|\dots|A_k])+d_{int}(\sum_{\small{\mbox{$(k+1)$-cuts}}}[A_1|\dots|A_{k+1}])=0.$$
\qed

We will need the following lemma in section 4.
\begin{lem}\label{generator}
For $A=(a_0,a_1,\dots,a_n;a_{n+1})$ let $\cc_A$ be the sub DGA of $\cc_S$ generated by
$A_1,\dots,A_k$ for all cuts of $A$. Then $H^0B(\cc_A)$ is generated by $T(A_i)$'s as an algebra.
\end{lem}
{\bf Proof.} It can be proved via a counting argument. Details are left to the reader.\qed

\begin{prop}
The coproduct of $T(A)$ is given by:
$$1\otimes T(A)+T(A)\otimes 1+\sum T(A_1)\otimes T(A_2)\cdots T(A_k)$$
where the sum is over all elementary $k$-cuts of $A$ for $k>1$. An elementary $k$-cut $A_1,\dots,A_k$ is a cut with
the property that each $A_i$ for $i>1$ is formed of consecutive elements, i.e. for some $l\ge 0$ and $m>1$ it is of the form  $(a_l;a_{l+1},\dots,a_{l+m-1};a_{l+m})$.
\end{prop}
\begin{rem}
This formula is analogous to Goncharov's formula for the coproduct of the framed mixed Hodge-Tate structure associated to iterated integrals in \cite{G}, theorem 6.1.
\end{rem}
{\bf Proof.} By definition of coproduct
$$\Delta T(A)=1\otimes T(A)+ T(A)\otimes 1+\sum [B_1|\dots|B_m]\otimes [B_{m+1}|\dots|B_l]$$
where the sum is over all $l$-cuts of $A$ with $l>1$ and $0<m<l$. Fix an elementary $k$-cut $c=(A_1,\dots,A_k)$ of $A$. A term $[B_1|\dots|B_m]\otimes [B_{m+1}|\dots|B_l]$
is said to be of type $c$ is $B_1,\dots,B_m$ is an $m$-cut of $A_1$. This relation partitions the terms of the above sum into groups indexed by elementary cuts of $A$.  If $[B_1|\dots|B_m]\otimes[B_{m+1}|\dots|B_l]$ is of type $c=(A_1,\dots,A_k)$ then obviously $B_{m+1},\dots,B_l$ is a shuffle permutation of cuts of $A_2,\dots,A_k$ and therefore the sum of terms of type $c$ gives $T(A_1)\otimes T(A_2)\cdots T(A_k)$.\qed

Let $\tilde{\cc}_S$ be a quotient of $\cc_S$ by the following two types of relations:
\begin{enumerate}
\item $(a_0;a_1,\dots,a_n;a_{n+1})+(a_{n+1},a_1,\dots,a_n;a_{n+2})=(a_0;a_1,\dots,a_n;a_{n+2}).$
\item $\sum_{\sigma\in\Sigma_{n,m}}(a_0;a_{\sigma(1)},\dots,a_{\sigma(n+m)};a_{n+m+1})=0.$
\item $(a_0;a_1,\dots,a_n;a_{n+1})=(-1)^{n}(a_{n+1};a_n,\dots,a_1;a_0).$
\end{enumerate}
Here $\Sigma_{n,m}$ is the subgroup of $(n,m)$ shuffles in $\Sigma_{n+m}$ consisting of those elements $\sigma$ such that:
$\sigma^{-1}(1)<\dots<\sigma^{-1}(n)$ and  $\sigma^{-1}(n+1)<\dots<\sigma^{-1}(n+m)$.
The ideal generated by these elements is both homogeneous and is stable under the differential and therefore $\tilde{\cc_S}$ is a DGA. Let $\tilde{T}(A)\in B(\tilde{\cc}_S)$ be the image of $T(A)$.
\begin{prop} \label{shuffle}We have the following relations for $\tilde{T}(A)$:
\begin{align*}
\sum_{i=0}^n \tilde{T}(a_0;a_1,\dots,a_i;a_{n+1})\tilde{T}(a_{n+1};a_{i+1},\dots,a_n;a_{n+2})&=\tilde{T}(a_0;a_1,\dots,a_n;a_{n+2})
\\
\tilde{T}(a_0;a_1,\dots,a_n;a_{n+m+1})\tilde{T}(a_0;a_{n+1},\dots,a_{n+m};a_{n+m+1})&=
\\
\sum_{\sigma\in \Sigma_{n,m}}
\tilde{T}(a_0;a_{\sigma(1)},&\dots,a_{\sigma(n+m)};a_{n+m+1})
\\
\tilde{T}(a_0;a_1,\dots,a_n;a_{n+1})= (-1)^n \tilde{T}(a_{n+1};a_n,\dots,a_1;a_0)
\end{align*}
For the first relation, we let $T(a_0;a_{n+1})=T(a_{n+1};a_{n+2})=1$.
\end{prop}
\begin{rem}
These formulae are analogous to the path composition and shuffle product relation for iterated integrals.
\end{rem}
{\bf Proof.} To prove the first relation, let $A=(a_0;\dots;a_{n+2}), B^i=(a_0;a_1,\dots,a_i;a_{n+1})$ and $C^i=(a_{n+1};a_{i+1},\dots,a_n;a_{n+2})$. Let $A_1,\dots,A_m$ be a cut of $A$ such that $A_1=A_{i_1},\dots, A_{i_k}$ end with $a_{n+2}$. One can replace $A_{i_j}=(a_l;\dots;a_{n+2})$ by the sum $A_{i_j}'+A_{i_j}''$ where $A_{i_j}'=(a_l;\dots;a_{n+1})$ and $A''_{i_j}=(a_{n+1};\dots;a_{n+2})$. Therefore $[A_1|\dots|A_k]$ can be written as:
\[[A_{i_1}'|\dots|\dots|A_{i_{k-1}}'|\dots|A_{i_k}'|A_{i_k+1}|\dots|A_m]+\]
\[\sum_{l=1}^k [A_{i_1}'|\dots|\dots|A_{i_{l-1}}'|\dots|A_{i_l}''|A_{i_l+1}|\dots|A_m]\]
If $A_{i_l}''=(a_{n+1};a_{i+1},\dots;a_{n+2})$, this element belongs to a term in $T(B^i)T(C^i)$. This way we get a 1-1 correspondence between the terms in $T(A)$ and the terms in $\sum_{i=0}^n T(B^i)T(C^i)$.

To prove the second relation, observe that the right hand side can be written as:
$$\sum_k [A_1|A_2|\dots|A_k]$$
where $(A_1,\dots,A_k)$ runs through all $k$-cuts of $A=(a_0;a_{\sigma(1)},\dots,a_{\sigma(n+m)},a_{n+m+1})$ for some $\sigma\in\Sigma_{n,m}$. If $A_i=(a_{i_0};a_{i_1},\dots,a_{i_l};a_{i_{l+1}})$ we let the support of $A_i$ be the set $\{i_1,\dots,i_l\}$. We say $A_i$ is mixed if its support has both elements $\le n$ and $>n$. We say $A_i$ is sequential if there is a $\sigma\in \Sigma_{n,m}$ and a subsequence $A_i'$, such that $(A_i',A_i)$ (in this order) is a 2-cut of $(a_0;a_{\sigma(1)},\dots,a_{\sigma(n+m)};a_{n+m+1})$.  Note that the last piece $A_k$ of a $k$-cut of $(a_0;a_{\sigma(1)},\dots,a_{\sigma(n+m)};a_{n+m+1})$ is always sequential. If $A_k$ is mixed and we sum over all shuffles of its support $\{i_1,\dots,i_l\}$ of elements $\le n$ and $>n$, leaving all other indices fixed, we find out that the terms $[A_1|\dots|A_k]$ corresponding to this cut, will cancel. So $A_k$ is pure, i.e.
$$A_k=(a_{p};a_i,a_{i+1},\dots,a_{i+j-1};a_q)$$
where $(i,i+1,\dots,i+j-1)$  are either all $\le n$ or all $>n$. Assume they are all $\le n$. I claim that we can also assume without loss of generality that $p=i-1$ and $q=i+j$ if $i+j\le n$ or $q=m+n+1$ if $i+j=n+1$, without changing the other pieces $A_1,\dots,A_{k-1}$. We give the arguments for $p$, the other index has a similar reasoning. If $p\ne i-1$ then it has to be bigger than $n$ so we can swap $a_p$ with $a_i,\dots,a_{i+j-1}$. Now the index before $p$ is either ${i-1}$ or ${p-1}$. In the first case by removing $a_i,\dots,a_{i+j-1}$ and doing the same cuts as before from the new sequence we get:
$$A_1,\dots,A_{k-1},(a_{i-1};a_i,\dots,a_{i+j-1};a_p)$$
But by the first type relations we have:
$$(a_{i-1};a_i,\dots,a_{i+j-1};a_p)+(a_p,a_i\dots,a_{i+j-1};a_q)=(a_{i-1};a_i,\dots,a_{i+j-1};a_q)$$
and so we have changed $p$ to $i-1$ without changing the rest of the cut. In the second case we can replace $A_k$ with $(a_{p-1};a_i,\dots,a_{i+j-1};a_q)$ and we can repeat the argument by decreasing $p$ until we get to $i-1$. Now we pick $A_{k-1}$; a similar argument shows that if its support is mixed then after summing over all the shuffles fixing the rest of the indices it will be canceled. And if it is pure, we can assume that its two ends are also in the same range. Continuing this way we can replace the terms in the right hand side with the terms in the left hand side. A neater argument can be given by an induction on $n+m$.
\\
The last relation is easy, and its proof is omitted.
\qed
\begin{rem}
It is interesting to observe that for proving the path composition we only used the first type relations and naively one expects that for the shuffle relation we only need the second type relations. However our proof uses both relations in order to prove the shuffle relation.
\end{rem}

We need the following lemma to deal with divergent iterated integrals
$$\int_{a_0}^{a_{n+1}}\frac{dt}{t-a_1}\circ \cdots\circ \frac{dt}{t-a_n}$$
where $a_0=a_1$ or $a_n=a_{n+1}$.
\begin{lem}\label{convergent}
The image of $\tilde{T}:\tilde{\cc}_S\to B(\cc_S)^0$ is generated by $\tilde{T}(a_0;a_1,\dots,a_n;a_{n+1})$ with $a_0\ne a_1$ and $a_n\ne a_{n+1}$ (by analogy we call these elements convergent) and $\tilde{T}(b_0;b_0;b_1)$.
\end{lem}
{\bf Proof.} Note that the shuffle relation implies:
$$\tilde{T}(a_0;\underbrace{a_1,a_1,\dots,a_1}_n;a_{n+1})=\frac{1}{n!}(\tilde{T}(a_0;a_1;a_{n+1}))^n$$
We prove the lemma for the following type elements by an induction on $i+j$:
$$\tilde{T}(a_0;\underbrace{a_0,\dots,a_0}_i,a_{i+1},\dots,a_{n-j},\underbrace{a_{n+1},\dots,a_{n+1}}_j;a_{n+1})$$
where $a_{i+1}\ne a_0$ and $a_{n-j}\ne a_{n+1}$, we may also assume $a_0\ne a_{n+1}$. Consider the product:
$$\tilde{T}(a_0;\underbrace{a_0,\dots,a_0}_i;a_{n+1})\tilde{T}(a_0;a_{i+1},\dots,a_{n-j};a_{n+1})\tilde{T}(a_0;\underbrace{a_{n+1},\dots,a_{n+1}}_j;a_{n+1})$$
By shuffle relation it is a linear sum of elements of the form $\tilde{T}(A)$, including the above term. Now notice that other than that term all the other terms by induction hypothesis can be written in terms of the convergent elements. \qed
\\

Let us define an ind-pro-object $\Pi^C(S,a,b)$ for $a,b\in S'$ with $S\subseteq S'$, in the category of graded $\h^0B({\cc}_{S,S'})$ comodule. Here $\cc_{S,S'}$ is the sub algebra of $\cc_{S'}$ generated by
$$(a_0;a_1,\dots,a_n;a_{n+1})\quad\mbox{such that}\quad  a_i\in S\quad i=1,\dots,n.$$
When $S$ is a finite subset of a field $F$, this will be a combinatorial analogue of the motivic torsor of paths $\Pi^{uni}({\Bbb A}^1\backslash S,a,b)$.
\\
A graded $\h^0B({\cc}_{S,S'})$ comodule is a graded $\Q$-vector space $V$ and a graded linear map $\nu : V\to V\otimes \h^0B(\cc_{S,S'})$ such that:
$$(id\otimes \Delta)\circ \nu = (\nu \otimes id)\circ \nu$$
$$(id\otimes \epsilon)\circ \nu =id$$
\begin{defn}\label{comb}
The ind-pro-object $\Pi^C(S,a,b)$ is the graded $\Q$-vector spaces of non-commuting power series with variables $X_s$, one for each element of $S$, $\Q\langle\langle X_s\rangle\rangle_{s\in S}$. This is an ind-object with respect to $S$ and a pro-object with respect to the polynomial degree or Adams grading. The grading is obtained by giving each $X_s$ weight $2$. For a sequence $A=(a;a_1,\dots,a_n;b)\in\cc_{S,S'}$ let $X_A=X_{a_1}\dots X_{a_n}$. Define:
 $$\nu: \Pi^C(S,a,b)\To \Pi^C(S,a,b)\otimes B({\cc}_{S,S'})^0$$
$$X_A\mapsto \sum X_{A_1}\otimes [A_2|\cdots|A_k]$$
where we sum over all cuts of $A$
\end{defn}
If we let the graded Adams piece $\h^0B(\cc_{S,S'})(r)$ to have weight $2r$, it is readily checked that above map is a graded linear map.
\begin{prop} $\Pi^C(S,a,b)$ is an ind-pro-object in the category of graded $H^0B(\cc_{S,S'})$ comodules. In particular $\nu$ factors through $\Pi^C(S,a,b)\otimes H^0B(\cc_{S,S'})$. If we make the construction with $\tilde{\cc}_S$, instead of $\cc_S$, there is a natural (path composition) morphism:
$$\Pi^C(S,a,b)\otimes \Pi^C(S,b,c)\To \Pi^C(S,a,c)$$
in the category of ind-pro graded $H^0B(\tilde{\cc}_{S,S'})$ comodules.
\end{prop}
{\bf Proof.} The coaction can be written as:
$$X_A\mapsto \sum X_{A_1}\otimes T(A_2)\cdots T(A_k)$$
where we sum over all elementary cuts of $A$. The rest follows directly from coproduct formula and path composition.\qed

\section{The Cycle map}

Let $F$ be a field. Our goal is to construct a natural map of DGA's:
$$\rho :{\cc}_F \To \N_F.$$
This induces a map
$$\rho: \h^0B({\cc}_F)\To \h^0B(\N_F)$$

\begin{defn} For $A=(a_0;a_1,\dots,a_n;a_{n+1})\in \cc_F$, relative to $\rho$ we define the motivic integral  $\I_{\rho}(A)\in \h^0B(\N_F)$ by
$\rho(\tilde{T}(A)).$
\end{defn}
\begin{rem} Since $\int_a^a\omega=0$, we expect to have $\rho(a_0;a_1,\dots,a_n;a_0)=0$. Therefore we will replace
$\cc_F$ with its quotient by the ideal generated by $(a_0;a_1,\dots,a_n;a_0)$.
\end{rem}
We confine ourself to construct $\rho$ on one of the following three sub-DGA's of $\cc_F$:
\begin{enumerate}
\item The sub-DGA $\cc_F^a$ of $a$-generic elements is generated by:
 $$(a_0;a_1,\dots,a_n;a_{n+1})\quad\mbox{such that}\quad a_i=a_j\:\:\mbox{and} \quad i\ne j \Rightarrow a_i=a.$$
\item The sub-DGA $\cc_F^1$ of sequentially distinct elements\footnote{Note that if $A$ is sequentially distinct and $A',A''$ is a 2-cut, then $A'$ might not be sequentially distinct, however in this case $A''=(a;\dots;a)$, an element that is set to be zero, so $\cc_F^1$ is really a sub-DGA.} is generated by:
$$(a_0;a_1,\dots,a_n;a_{n+1})\quad\mbox{such that}\quad a_i\ne a_{i+1}\quad\mbox{for}\quad i=1,\dots,n.$$
\item The sub-DGA $\cc_F^2$ of binary elements is generated by:
$$(a_0;a_1,\dots,a_n;a_{n+1})\quad\mbox{such that}\quad a_i\in\{0,1\}\quad \mbox{for}\quad i=1,\dots,n.$$
\end{enumerate}
The technical definition that follows, is needed for the Hodge realization. We want to construct cycles over the field
$F(s_1,\dots,s_k)$. For fixed base points $a_0$ and $a_{n+1}$, there is a (partially defined) map:
$$\delta:\N_{F(s_1,\dots,s_k)}\To \N_{F(s_1,\dots,s_{k-1})}$$
\begin{align*}
(\delta Z)(s_1,\dots,s_{k-1})&=Z(a_0,s_1,\dots,s_{k-1})-Z(s_1,s_1,\dots,s_{k-1})+\dots
\\
+(-1)^{k-1} &Z(s_1,\dots,s_{k-2},s_{k-1},s_{k-1})+(-1)^kZ(s_1,\dots,s_{k-1},a_{n+1})
\end{align*}
 We say partially defined because the terms of this sum should be admissible. The cycles that we construct have very simple (in fact linear) coefficients in terms of $s_i$'s and so these specializations are sensible.
\begin{defn}\label{int}
An integration theory on a sub-DGA $\cc_F'$ of $\cc_F$, is a sequence of algebra morphism (not DGA) for $k\ge 1$:
$$\rho_k:\cc_F'\To \N_{F(s_1,\dots,s_k)}[k-1]$$ where the bracket denotes a shift in the cohomological degree.
i.e. if $A\in \cc_F^{'1}(n)$ then $\rho_k(A)\in \N_{F(s_1,\dots,s_k)}^k(n)$,  is a codimension $n$ cycle inside $\Box^{2n-k}$. Let $A=(a_0;a_1,\dots,a_n;a_{n+1})\in \cc_F'$. The following properties should hold:
\begin{enumerate}
\item $\rho_k(A)=0$ for $k>n$.
\item $\rho_n(A)=\mbox{Alt}(c_1({s_1}-{a_1}),\dots,c_n(s_n-a_n))$\footnote{All is needed for the next section is: $\int_{\eta_1}\frac{dz_1}{z_1}\wedge\dots\wedge\frac{dz_n}{z_n}= \int_{\gamma}\frac{dt}{t-a_1}\circ\dots\circ \frac{dt}{t-a_n}$ where $\eta_1$ is the topological cycle obtained from $\rho_1(A)$ by replacing $s_i$'s with $\gamma(s_i)$'s with $0\le s_1\dots\le s_n\le 1$ for a path $\gamma$ from $a_0$ to $a_{n+1}$ in $\C-\{a_1,\dots,a_n\}$.}
for some non-zero constants $c_i$'s that might depend on $A$.
\item The differential of $\rho_k$ is: $$d\rho_k(A)= -\delta\rho_{k+1}(A)+(-1)^k\sum_{\mbox{\small 2-cuts of A}}\rho_k(A')\delta\rho_1(A'')$$
Implicit in this definition is admissibility of the cycles $\rho_k(A)$ and $\delta\rho_k(A)$ for $(k,n)\ne (1,1)$.
The case $\delta\rho_1(a_0;a_1;a_2)$ needs special definition when $a_0=a_1$ or $a_1=a_2$.
\end{enumerate}
\end{defn}
\begin{lem}\label{lem2}
The differential of $\delta\rho_k(A)$ is given by:
 $$d(\delta\rho_k(A)) = (-1)^{k}\sum_{\small{\mbox  2-cuts of A}}   \delta\rho_k(A')\delta\rho_1(A'').$$               \end{lem}                                                                                                            {
{\bf Proof.} Take boundary of both sides of equation in condition (3). Use the evident facts that $\delta d=d\delta$
and $\delta^2=0$ and $\delta(a\delta b)=(\delta a)(\delta b)$.\qed
\begin{rem}
According to lemma \ref{lem2} the map $\rho := \delta\rho_1 :\cc_F'\To \N$ is a morphism of DGA's. So if we have an integration theory we have an associated motivic integral in $\h^0B(\N)$. In the next section we show that it has the right Hodge realization, regardless of the integration theory used.
\end{rem}
\begin{defn}
An integration theory is called base point free, if $\rho_k(a_0;\dots;a_{n+1})$ is independent of $a_0$ and $a_{n+1}$. So for such theories we sometimes use the notation $\rho_k(a_1,\dots,a_n)$. An integration theory is called permuting if it is base point free and for all $n,m>0$:
$$\sum_{\sigma\in\Sigma_{n,m}}\rho_1(a_{\sigma(1)},\dots,a_{\sigma(n+m)})=0$$
$$\rho_1(a_1, a_2,\dots,a_n)=(-1)^{n-1}\rho_1(a_n,a_{n-1},\dots,a_1)$$
If we have a permuting theory the corresponding $\rho:\cc_F\To\N$ factors through $\tilde{\cc}_F$, and hence its motivic integrals $\I(A)$ satisfy the shuffle, path composition and inversion relations similar to those in proposition \ref{shuffle}.
\end{defn}

The construction of an integration theory is achieved inductively.

$$\hat{\rho}_1(a_1)(s_1)=\begin{cases} \mbox{Alt}(1-\frac{s_1}{a_1})&\quad\mbox{if}\quad a_1\ne 0
\\
\mbox{Alt}(s_1)&\quad\mbox{if}\quad a_1=0
\end{cases}
$$
 We also let :
$$Sp|_{s_1=a_0}\hat{\rho}_1(a_1)(s_1)=\begin{cases} \hat{\rho}_1(a_1)(a_0) &\quad\mbox{if}\quad a_0\ne a_1
\\
\hat{\rho}_1(a_1)(1+a_0)&\quad\mbox{if}\quad a_0=a_1
\end{cases}
$$
This specialization is used to regularize divergent integrals. They reflect the fact that with the standard tangential base point at $0$ one should define $\int_0^1\frac{dz}{z}$ to be zero. For any integration theory the most important part is a definition for $\rho_1(A)$, for $A=(a_0;a_1,\dots,a_n;a_{n+1})$.
In fact if we let $A_i=(a_0;a_1,\dots,a_i;a_{i+1})$ and $A_i'=(a_{i};a_{i+1},\dots,a_n;a_{n+1})$, we can define $\rho_k$ inductively:
$$\rho_k(A)(s_1,\dots,s_k)=\sum_{i=1}^{n-k+1}\rho_1(A_i)(s_1)\cdot \rho_{k-1}(A_i')(s_2,\dots,s_k)$$
\begin{lem}\label{easy}
If $\rho_1(A)$ satisfies $d\rho_1(A)=-\delta\rho_2(A)-\sum\rho_1(A')\delta\rho_1(A'')$, with $\rho_k$'s defined as above in terms of $\rho_1$ for all $A$ of depth $<n$. Then we have
$$d\rho_k(A)=-\delta\rho_{k+1}(A)+(-1)^k\sum \rho_k(A')\delta\rho_1(A'')$$
for all $A$ of depth $<n+k-1$.
\end{lem}
{\bf Proof.} By definition:
$$d\rho_k(A)=\sum_{i=1}^{m-k+1} d\rho_1(A_i)\cdot \rho_{k-1}(A_i')-\rho_1(A_i)\cdot d\rho_{k-1}(A_i')$$
where $m$ is the depth of $A$, if $m < n+k-1$ then $i<n$ so we can use:
$$d\rho_1(A_i)=-\delta\rho_2(A_i)-\sum_{\mbox{\tiny{2-cuts of $A_i$}}} \rho_1(A_{i1})\delta\rho_1(A_{i2})$$
and inductively:
$$d\rho_{k-1}(A_i')=-\delta\rho_k(A_i')-\sum_{\mbox{\tiny{2-cuts of $A_i'$}}} \rho_{k-1}(A_{i1}')\delta\rho_1(A_{i2}')$$
so $d\rho_k(A)$ is the sum of the following two terms:
$$\sum -\delta\rho_2(A_i)\cdot \rho_{k-1}(A_i')+\rho_1(A_i)\cdot\delta\rho_k(A_i')$$
$$(-1)^k\sum_i\left(\sum \rho_1(A_{i1})\rho_{k-1}(A_i')\delta\rho_1(A_{i2}) +\sum \rho_1(A_i)\rho_{k-1}(A_{i1}')\delta_1\rho(A_{i2}')\right)$$
One can check that by definition of $\rho_k$ and fixing $A_{i2}$ or $A_{i2}'$ and letting $i$ vary that the second sum is just
$$(-1)^k\sum_{\mbox{\tiny 2-cuts of A}}\rho_k(A')\delta\rho_1(A'')$$
For the first sum if we neglect the two terms $\rho_2(A_i)(s,a_{i+1})$ and $\rho_{k}(A_i')(a_i,s_1,\dots,s_{k-1})$
we get $-\delta\rho_{k+1}(A)$, now one sees immediately that these two terms when summed over $i$ cancel each other.\qed

Let
$$\rho_1(a_0;a_1;a_2)(s)=\begin{cases}
\hat{\rho}_1(a_1-a)(s-a)&\mbox{for $a$-generic case}
\\
\hat{\rho}_1(a_1-a_0)(s-a_0)&\mbox{for sequentially distinct case}
\\
\hat{\rho}_1(a_1)(s)&\mbox{for binary case}
\end{cases}
$$
Note that the first and last case are base point independent. These come with their specialization coming from $\hat{\rho}_1$. We set:
$$\rho(a_0;a_1;a_2):= \delta\rho_1:= Sp|_{s=a_0}\rho_1(a_0;a_1;a_2)(s)-Sp|_{s=a_2}\rho_1(a_0;a_1;a_2)(s).$$
To define $\rho_1(A)$ for $n>1$, we assume that we have defined $\rho_1$ for all elements of depth $<n$ so we have a definition of $\rho_k$ for all elements of depth $<n+k-1$. For the $a$-generic and the sequentially distinct case, let:
$$\rho_1(A)(s)=\mbox{Alt}\left(\frac{s-t}{c-t},\rho_2(t,t)\right)$$
where $c=a$ in the $a$-generic case and $c=a_0$ for the sequentially distinct case.
Here we let $t$ vary. In fact inductively one sees that $\rho_1(A)$ has $n-1$ of these internal variables, i.e., the way $\rho_1$ is defined gives us a (formal) linear combination of rational functions over $F(s)$: $$\phi_i: ({\mathbb P}^1)^{n-1}\To ({\mathbb P}^1)^{2n-1}.$$
The cycle is the locus of these functions, intersected with $\Box^{2n-1}$. More precisely we take the pull-back of these function under the inclusion $\Box^{2n-1}\hookrightarrow ({\mathbb P}^1)^{2n-1}$ to get
$$\phi_i':X_i\To \Box^{2n-1}\quad\quad X_i=({\mathbb P}^1)^{n-1}\times_{({\mathbb P}^1)^{2n-1}} \Box^{2n-1}$$
These maps are proper so we have push forward $\phi_{i*}'$ for cycles and the cycle is:
$$\rho_1(A) = \sum_i \phi_{i*}'(X_i).$$

The definition for the binary case is more involved. One sees immediately that $\rho_1$ is base point dependent (in fact just the initial base point) for the sequentially distinct case and it is base point independent for the $a$-generic case.

\begin{prop}
The morphisms $\rho_k$  define an integration theory for $a$-generic and sequentially distinct sub algebras . For $a$-generic the theory is permuting, hence we have shuffle, path composition and inversion relations for $\I_{\rho}$.
\end{prop}
{\bf Proof.} Assume $\rho_1(A)$ is admissible and satisfies
$$d\rho_1(A)=-\delta\rho_2(A)-\sum\rho_1(A')\delta\rho_1(A'').$$
 for $A$ of depth $<n$. By lemma \ref{easy} for $A$ of depth $n$ we have:
$$d\rho_2(A)=-\delta\rho_3(A)+\sum\rho_2(A')\delta\rho_1(A'')$$
but if we specialize at $(t,t)$ we can replace
$-\delta\rho_3(t,t)$ with $-\rho_3(a_0,t,t)+\rho_3(t,t,a_{n+1})$
Assuming the admissibility, differential of $\rho_1(A)=(\frac{s-t}{c-t},\rho_2(A)(t,t))$ is:
$$\rho_2(s,s)+(\frac{s-t}{c-t},\rho_3(a_0,t,t))-(\frac{s-t}{c-t},\rho_3(t,t,a_{n+1}))-(\frac{s-t}{c-t},\sum\rho_2(A')(t,t)\delta\rho_1(A''))$$
The first term is by $t=s$ for the first coordinate, also when $t=c$ the second factor becomes empty, This is because $\rho_2(A)(c,c)$ has at least one coordinate which is $1$. The change in sign happens since the coordinates of $\rho_2$ are shifted by 1 in the definition of $\rho_1$. Now if we use the definitions we see that this is exactly the formula we wanted to prove for $d\rho_1$. It is left to check that $\rho_1(A)$ is admissible, i.e. its intersection with any face $\mathcal F$ of $\Box^{2n-1}$ is either empty or is of codimension $n$. First let's prove inductively that the intersection of $\rho_1\cap {\mathcal F}$ is empty of we set $s=c$ ($c$ is $a_0$ for sequentially distinct and is $a$ for the generic). Recall that $\rho_1(A)(s)=\mbox{Alt}(\frac{s-t}{c-t},\rho_2(A)(t,t))$. It is enough to prove this claim before alteration. If $\mathcal F$ contains $z_1=\infty$ then we have to set $t=c$ in $\rho_2(A)(t,t)$ which by induction will give empty. If the face doesn't contain $z_1=0$ then $\frac{s-t}{c-t}$ survives after intersection and if we let $s=c$ it is 1 so the intersection with $\Box^m$ is empty. Finally if it does have $z_1=0$ we have to take the intersection $\rho_2(s,s)$ which is $\sum \rho_1(A_i)(s)\rho_1(A_i')(s)$ and by induction this becomes empty if we set $s=c$. Next we show that the intersection of $\rho_1$ with any face $z_i=\infty$ is empty. We did the case $z_1=\infty$ before. For other cases we have to take the intersection of $\rho_2(t,t)$ with an infinite face, if this is non-empty by induction we need to have $t=\infty$ which makes $\frac{s-t}{c-t}$ to be 1. The best way to prove the admissibility in general is to use the language of trees. We refer the reader to \cite{FJ} or \cite{GGL} for this. For a rooted tree $T$ (not necessarily binary) with leaves labeled with $a_1,\dots,a_n$ and $E$ edges, define a subvariety of $\Box^{E}$ by:
$$\rho_1^T(A)(s)=(\frac{s-t}{c-t},\rho_1^{T_1}(supp(T_1))(t),\dots,\rho_1^{T_k}(supp(T_k))(t))$$
where $T_1,T_2,\dots,T_k$ are the children sub-trees of the root from left to right. And if $T_i$ has leaves labeled with $a_j,a_{j+1},\dots,a_m$ then $supp(T_i)=(a_{j-1};a_j,\dots,a_m;a_{m+1})$. It is observed readily that:
$$\rho_1(A)(s)=\sum \rho_1^T(A)(s)$$
where we sum over the binary trees with leaves labeled at $a_1,\dots,a_n$. Now note that the intersection of $\rho^T(A)$ with $z_i=0$ where the coordinate corresponds to an internal edge $e$ of $T$ gives $\rho_1^{T'}(s)$ where $T'$ is obtained from $T$ by contracting $e$. And its intersection with
$z_i=0$ corresponding to a leaf labeled at $a_i$ and with its sibling sub-trees $T_1,\dots,T_k$ corresponds to
$$\pm \rho_1^{T'}(supp(T'))(s)\rho^{T_1}(supp(T_1))(a_i)\dots\rho^{T_k}_1(supp(T_k))(a_i).$$
where $T'$ is obtained from $T$ by contracting the leaf with label $a_i$ and cutting off all its siblings and labeling the new made leaf by $a_i$. So if we want to use induction on the number of edges, we only need to show that the case where one of the $T_i$'s above is a leaf with label at $a_i$ will give an empty cycle. In the case of $a$-generic we need to have $a_i=a$ and in this case at least one other $T_i$ gives $1$ if we replace its $s$ variable with $a$ (note that if all $T_i$'s are leaves with labels at $a$ the alteration will kill this term.) for the sequentially distinct case, the sibling $T_1$ to the right of the contracted leaf at $a_i$ has its support of the form $(a_i;a_{i+1},\dots;a_{m+1})$ and hence $\rho_1^{T_1}(supp(T_1))(a_i)$ is empty inside $\Box^l$. This finishes the admissibility condition.
\\
The fact that the integration theory on $a$-generic elements is permuting is left to the reader.
\qed

Now we deal with the binary case. We have to define $\rho_1(a_0;a_1,\dots,a_n;a_{n+1})$, but since the theory that we will construct is base point independent, we write $\rho_1(A)$ for $A=(a_1,\dots,a_n)$. We need axillary cycles $\rho_1^{\epsilon}(A)$ and $\rho_1^{\epsilon}(A_i|A_i')$ for ${\epsilon} = 0$ or $1$ and $A_i=(a_1,\dots,a_i)$, $A_i'=(a_{i+1},\dots,a_n)$. The definitions are inductive. First we give the special and simple case of $n=2$:
$$\rho_1^{\epsilon}(A)=\rho_1^{\epsilon}(A_1|A_1')=\frac{1}{2}\left(\frac{s-t}{\epsilon-t},\rho_1(a_1)(t),\rho_1(a_2)(t)\right)$$
Note that by alteration this term is zero when $a_1=a_2$.
To make our formulation more compact, we also set:
$$\rho_1^{\epsilon}(a)(s)=|{\epsilon}-a|\rho_1(a)(s)=\begin{cases}
0&\mbox{if}\quad a=\epsilon
\\
\mbox{Alt}(s)&\mbox{if}\quad 0=a\ne \epsilon
\\
\mbox{Alt}(1-s) &\mbox{if}\quad 1=a\ne\epsilon
\end{cases}
$$
Remind that 0 is the empty cycle and not zero as a function. Now the case $n>2$:
$$\rho_1^{\epsilon}(A)=\mbox{Alt}\sum_{i=1}^{n-1}\rho_1^{\epsilon}(A_i|A_i')-\rho_1^{\epsilon}(A_i'|A_i)$$
\begin{align*}
\rho_1^{\epsilon}(A_i|A_i')=&\left(\frac{s-t}{{\epsilon}-t},\rho_1^{1}(A_i)(t),\rho_1^{0}(A_i')(t)\right)+\delta_i\left(\frac{t-s}{\epsilon-t},\rho^{1}_1(A_i)(t),\rho^{0}_1(A_i')(s)\right)
\\
&+\delta_i\left(\frac{t-s}{\epsilon-t},\rho^{0}_1(A_i)(t),\rho^{1}_1(A_i')(s)\right)
\end{align*}
where $\delta_i=-\frac{1}{2}$ if $i\ne 1$ and $n-1$ and $\delta_1=0$, $\delta_{n-1}=-1$ (remind that $n>2$).
Now for $A=(a_1,\dots,a_n)$ let:
$$\rho_1(A)(s)=\rho_1^0(A)(s)+\rho_1^1(A)(s)$$
$$\rho_k(A)(s_1,\dots,s_k)=\sum_{i=1}^{n-k+1}\rho_1(A_i)(s_1)\rho_{k-1}(A_i')(s_2,\dots,s_k).$$
\begin{prop}
The morphisms $\rho_k$ defined above, form a permuting integration theory on the binary subalgebra $\cc_F^2$, generated by $(a_0;a_1,\dots,a_n;a_{n+1})$ with $a_i\in\{0,1\}$ for $i=1,\dots,n$.
\end{prop}
\begin{rem}
One can easily modify this construction to work for binary elements $(a_0;\dots;a_{n+1})$ with $a_i\in\{\alpha_0,\alpha_1\}$ for two fixed distinct elements $\alpha_0,\alpha_1\in F$. The starting point would be:
$$\rho_1(a)(s)=\begin{cases}
\frac{s-a}{\alpha_0-a}&\mbox{if}\quad a\ne \alpha_0
\\
\frac{s-a}{\alpha_1-a}&\mbox{if}\quad a=\alpha_0
\end{cases}
$$
Now similarly we define $\rho_1^{0}$ and $\rho^{1}_1$. Where instead of ${\epsilon}$, we have $\alpha_{\epsilon}$.
\end{rem}
{\bf Proof.} According to lemma \ref{easy} all we need to do is to check $\rho_1$ is admissible and satisfies
$$d\rho_1=-\delta\rho_2-\sum_{\mbox{\tiny{2-cuts of A}}}\rho_1(A')\delta\rho_1(A'')$$
The admissibility follows from the following observation:
\\
The intersection ${\mathcal F}\cap \rho_1^{\epsilon}(A)(s)$ of $\rho_1^{\epsilon}$ with any face of $\Box^{2n-1}$ is empty if we set $s=\epsilon$. Note that this is stronger than just saying $\rho_1^{\epsilon}(\epsilon)$ is empty, which in fact is obvious since the first coordinate becomes 1. this stronger version follows since each component has both $\rho_1^0$ and $\rho_1^1$, if we intersect with $z_1=0$ to kill the obvious coordinate giving $1$ we still have persisting coordinates that gives $1$, now regardless of $\epsilon$. If we intersect with $z_1=\infty$ we get the two terms $\rho_1^1(A_i)(\epsilon)\rho_1^0(A_i')(s)$ and $\rho_1^0(A_i)(\epsilon)\rho_1^1(A_i')(s)$. Now if $\epsilon=0$ the second component is empty by induction and now for the first component we can use inductive hypothesis to show any intersection with a face becomes empty if we set $s=0$, the case $\epsilon=1$ is similar.
This will imply the admissibility.

Now again as before the definition can be restated in terms of rooted binary trees with leaves at $a_1,\dots,a_n$. Such a tree gives a partition of this set into $A_i$ and $A_i'$ (left and right child labels). And to define $\rho_1$
or $\rho_1^{\epsilon}$ we need this partition. In fact for a binary tree $T$ with left and right subtrees $L$ and $R$ we define:
\begin{align*}
\rho_1^{\epsilon}(T)=&\left(\frac{s-t}{{\epsilon}-t},\rho_1^{1}(L)(t),\rho_1^{0}(R)(t)\right)+\delta_i\left(\frac{t-s}{\epsilon-t},\rho^{1}_1(L)(t),\rho^{0}_1(R)(s)\right)
\\
&+\delta_i\left(\frac{t-s}{\epsilon-t},\rho^{0}_1(L)(t),\rho^{1}_1(R)(s)\right)
+\left(\frac{s-t}{{\epsilon}-t},\rho_1^{0}(L)(t),\rho_1^{1}(R)(t)\right)
\\
&+\delta_{n-i}\left(\frac{t-s}{\epsilon-t},\rho^{0}_1(L)(s),\rho^{1}_1(R)(t)\right)
+\delta_{n-i}\left(\frac{t-s}{\epsilon-t},\rho^{1}_1(L)(s),\rho^{0}_1(R)(t)\right)
\end{align*}
with $\delta_i$ as before with $i$ being the number of leaves of $R$. With this definition we can write:
$$\rho_1^{\epsilon}(A)=\sum\rho_1^{\epsilon}(T)$$
where we sum over all binary trees with leaves labeled with $A=(a_1,\dots,a_n)$.
Now take an edge $e$ of one of these binary trees $T$, it corresponds to a coordinate $z_e$ of $\rho_1^{\epsilon}(A)$. If it is not a leaf, in 6 different ways i.e. one of the 6 terms above. If
one or both subtrees of this edge are leaves some of these 6 cases are empty. To manage these different choices we put an extra decoration on our binary trees. Namely we decorate each edge (excluding the leaves) with a number 1 through 6, corresponding to the above 6 terms in the order they are written. Now if $(T,dec_T)$ is such a decorated tree with root decorated by $1\le j\le 6$, we let $\rho_1^{\epsilon}(T,dec_T)$ be the $j$'th term of the above expression, where $L$ and $R$ get their decorations from $T$. We can rewrite the above formula this time with decorated binary trees:
$$\rho_1^{\epsilon}(A)=\sum\rho_1^{\epsilon}(T,dec_T)$$
The $0$-special edges of a decorated binary tree are the internal edge connected to the left of the root with the root decoration of $2$ or $3$ and the internal edge connected to the right of the root
with the root decoration of 5 or 6 and the root itself.
Similarly the $\infty$-special edges are the internal ones connected to a leaf with a decoration of $2,3,5$ or $6$.
\\
We want to show that there is an involution $(T, dec_T, e,\alpha)\leftrightarrow (T',dec_{T'}, e',\alpha')$ between quadruples a binary tree (with leaves at $a_1,\dots,a_n$) ,
its decoration, an edge $e$ and $\alpha=0$ or $\infty$ such that with the right sign (dictated by the differential on cycles)
$\rho_1^{\epsilon}(T,dec_T)|_{z_e=\alpha}$ cancels $\rho_1^{\epsilon}(T',dec_{T'})|_{z_{e'}=\alpha'}$, with the exclusion of $(T,dec_T,e,\alpha)$ where
$e$ is an $\alpha$-special edge. This is done by considering different cases.

 To calculate $d\rho^{\epsilon}(A)$ we only need to consider $\rho^{\epsilon}(T,dec_T)|_{z_e=\alpha}$ with $e$ being an $\alpha$-special edge for $\alpha=0$ or $\infty$.
So we have the following cases. For ease of notation we right $L^{\epsilon}$ for $\rho_1^{\epsilon}(L,dec_L)(s)$ and similarly for
$R$.
\begin{enumerate}
\item $e$ is the root with decoration 1: $$\rho_1^{\epsilon}(T,dec_T)|_{z_e=0}=L^1R^0.$$
\item $e$ is the root with decoration 4: $$\rho_1^{\epsilon}(T,dec_T)|_{z_e=0}=L^0R^1.$$
\item $e$ is the root with decoration 2 or 6 (note that $\delta_i+\delta_{n-i}=1$):
$$\rho_1^{\epsilon}(T,dec_T)|_{z_e=0}+\rho^{\epsilon}(T,dec_T')|_{z_e=0}=-L^1R^0.$$
where $dec_T$ and $dec_T'$ are identical except at the root, one with decoration 2 and one with decoration 6.
\item $e$ is the root with decoration 3 or 5:
$$\rho_1^{\epsilon}(T,dec_T)|_{z_e=0}+\rho_1^{\epsilon}(T,dec_T')|_{z_e=0}=-L^0R^1.$$
\item $e$ is to the left (or right) of the root with the root decoration 2 (or 6):
$$-\rho_1^{0}(T,dec_T)|_{z_e=0}+\rho_1^{0}(T,dec_T')|_{z_e=0}=L^0R^0.$$
\begin{align*}
-\rho_1^{1}(T,dec_T)|_{z_e=0}+\rho_1^{1}(T,dec_T')|_{z_e=0}=&-\delta_iL^1R^0-\delta_{n-i}L^0R^1
\end{align*}
\item $e$ is to the left (or right) of the root with the root decoration 3 (or 5):
$$-\rho_1^{1}(T,dec_T)|_{z_e=0}+\rho_1^{1}(T,dec_T')|_{z_e=0}=L^1R^1$$
\begin{align*}
-\rho_1^{0}(T,dec_T)|_{z_e=0}+\rho_1^{0}(T,dec_T')|_{z_e=0}=&-\delta_iL^0R^1-\delta_{n-i}L^1R^0
\end{align*}
\end{enumerate}
Adding all these terms will give us $L^0R^0+L^1R^0+L^1R^1+L^0R^1$ or $(L^0+L^1)(R^0+R^1)$ when we sum over all decorated trees
this gives $\rho_2(A)(s,s)$, which is a term in $-\delta\rho_2(A)$. The other terms in the formula:
$$d\rho_1(A)=-\delta\rho_2(A)+\sum\rho_1(A')\delta\rho_1(A'')$$
come from $\infty$-special edges.
\begin{enumerate}
\item The edge $e$ has a right child which is a leaf at $a_{i+1}$ and and a left child which is a tree $T'$ with labels
$A_i=(a_1,\dots,a_i)$, as shown below, with a decoration of 2 or 3 then:
$$\rho_1^{\epsilon}(T,dec_T)|_{z_e=\infty}+\rho_1^{\epsilon}(T,dec'_T)|_{z_e=\infty}=
-\rho_1(T',dec_{T'})(a_{i+1})\rho_1^{\epsilon}(T'',dec_{T''})(s)$$
where $T''$ is the tree obtained from $T$ by removing the leaf $a_{i+1}$ and its sibling subtree $T'$ and labeling the newly
made leaf by $a_{i+1}$.
\vskip 10pt \hskip 100pt
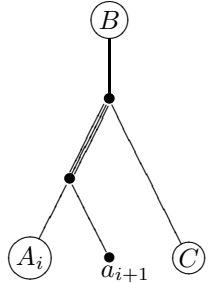
\hbox{
 \vbox{\xy 0;<-3pt,0pt>:
 \POS(0,0) *+{B} *\cir{}
 \ar @{-} +(0,10)*{\bullet}
  \POS(0,10) *{\bullet}
 \ar @{=} +(5,10)*{\bullet}
 \ar @{-} +(5,10)
 \ar @{-} +(-10,20)*+{C}*\cir{}
 \POS(5,20) *{\bullet}
 \ar @{-} +(5,10)*+{A_i}*\cir{}

 \ar @{-} +(-5,10)*{\bullet}
 \POS(-2,32) *{a_{i+1}}
 \endxy}
}
\vskip 10pt
If we add the trivial term
$$-\rho_1(T',dec_{T'})(a_0)\rho_1^{\epsilon}(T'',dec_{T''})(s)+\rho_1(T',dec_{T'})(a_0)\rho_1^{\epsilon}(T'',dec_{T''})(s)$$
and add over all possible cases of this form we get the term
$$-\rho_2(A)(a_0,s)-\sum_i\rho_1(a_{i+1},\dots,a_n)\delta\rho_1(a_0;a_1,\dots,a_i;a_{i+1})$$
\item The edge $e$ has a left child which is a leaf at $a_{i}$ and a right child which is a tree $T'$ with labels $A_i'=(a_{i+1},\dots,a_n)$
and decoration 4 or 5.
\vskip 10pt \hskip 100pt
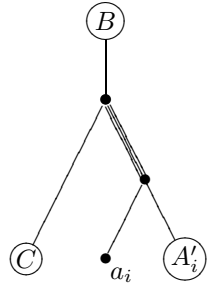
\hbox{
 \vbox{\xy 0;<-3pt,0pt>:
 \POS(0,0) *+{B} *\cir{}
 \ar @{-} +(0,10)*{\bullet}
  \POS(0,10) *{\bullet}
 \ar @{-} +(10,20)*+{C}*\cir{}
 \ar @{-} +(-10,20)*+{A_i'}*\cir{}
 \ar @{=} +(-5,10)
 \POS(-5,20)*{\bullet}
 \ar @{-} +(5,10)*{\bullet}
 \POS(-2,32) *{a_{i}}
 \endxy}
}
\vskip 10pt
A calculation similar to the above case gives:
$$-\rho_2(A)(s,a_{n+1})-\sum_i\rho_1(a_{1},\dots,a_i)\delta\rho_1(a_i;a_{i+1},\dots,a_n;a_{n+1})$$

\item One sees that the following two types of $\infty$-special edges with $A_{i,j}=(a_{i+1},\dots,a_j)$ will give the rest of the differential:
$$-\sum\rho_1(A')\delta\rho_1(A'')$$
where we sum over 2-cuts $A',A''$ of $(a_0;\dots;a_{n+1})$ such that the initial and final point of $A''$ are
$a_i$ and $a_{j+1}$ with $0<i<j<n$.

\vskip 10pt \hskip 100pt              
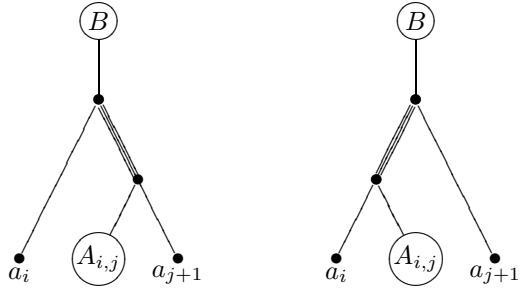
\hbox{
 \vbox{\xy 0;<-3pt,0pt>:
 \POS(0,0) *+{B} *\cir{}
 \ar @{-} +(0,10)*{\bullet}
  \POS(0,10) *{\bullet}
 \ar @{=} +(5,10)*{\bullet}
 \ar @{-} +(5,10)
 \ar @{-} +(-10,20)*{\bullet}
 \POS(5,20) *{\bullet}
 \ar @{-} +(5,10)*{\bullet}
 \ar @{-} +(-5,10)*+{A_{i,j}}*\cir{}
 \POS(10,32)*{a_i}
 \POS(-10,32)*{a_{j+1}}
 \POS(40,0) *+{B} *\cir{}
 \ar @{-} +(0,10)*{\bullet}
  \POS(40,10) *{\bullet}
 \ar @{-} +(10,20)*{\bullet}
 \ar @{-} +(-10,20)*{\bullet}
 \ar @{=} +(-5,10)
 \POS(35,20)*{\bullet}
 \ar @{-} +(5,10)*+{A_{i,j}}*\cir{}
 \POS(30,32) *{a_{j+1}}
 \POS(50,32) *{a_i}
 \endxy}
}
\vskip 10pt
\end{enumerate}
The proposition is proved.\qed
\\

In definition \ref{comb} we defined an ind-pro object $\Pi^C(S,a,b)$ for a given set $S\subseteq S'$ and a co-action:
$$\nu:\Pi^C(S,a,b)\To \Pi^C(S,a,b)\otimes \h^0B(\tilde{\cc}_{S,S'})$$
Now for $S=\{0,1\}$ and $S'=\{0,1,a,b\}$ we have a map from  $\tilde{\cc}_{S,S'}$ to $\N$ so we get a map
$$\Pi^C(S,a,b)\To \Pi^C(S,a,b)\otimes \h^0B(\N)$$
$$X_A\mapsto \sum X_{A_1}\otimes \I_{\rho}(A_2)\cdots \I_{\rho}(A_k)$$
where we sum over all elementary cuts of $A$.
Which makes $\Pi^C(S,a,b)=\Q\langle\langle X_0,X_1\rangle\rangle$ into a comodule over $\h^0B(\N)$ .
\begin{defn}
The motivic torsor of paths on ${\mathbb P}^1-\{0,1,\infty\}$ from $a$ to $b$ is the above co-module over $\h^0B(\N)$.
It is denoted by $\Pi^{Mot}({\mathbb P}^1-\{0,1,\infty\};a,b)$.
\end{defn}
In the next section we show that its Hodge realization is the canonical Hodge structure on the unipotent completion of the torsor of paths. (see appendix B for a review of it)

\section{Hodge Realization of the Motivic Iterated Integrals}

As in definition \ref{int}, let $\{\rho_k\}$ be an integration theory, and $\I(A)$ its relative motivic iterated integral.
In this section we show that the Hodge realization of $\I(A)$ is the canonical framed MHTS $\I^{\mathcal H}(A)$(refer to definition \ref{framed integral})
for iterated integrals.
\\
Given $A=(a_0; a_1,\dots, a_n; a_{n+1})$, let $\N'$ be the sub algebra of $\N$ generated by $\rho(A_i)$'s for all cuts
$A_1,\dots,A_k$ of $A$. We want to define a sub-DGA  $\D'$ of topological cycles that satisfies the condition of the definition \ref{admissible} in appendix B. To do this we need to define some intermediate topological cycles. Their definition depends on a choice of a path $\gamma$ with
interior in $\C-\{a_1,\dots,a_n\}$ from the tangential base point
$a_0$ to the tangential base point $a_{n+1}$.

 For $1\le i\le n$
define a topological cycles $\eta_i(A)$ inside
$\Box_{\C}^{2n-i}$ with (real) dimension $2n-i$ by:

$$\eta_i(A)=\rho_i(A)(\gamma(s_1),\dots,\gamma(s_i))$$
where $s_k$'s vary in the simplex $0\le s_1\le\dots\le s_i\le 1$.
\\
 Let $\Gamma$ be a small disk around zero in $\Box_{\C}$ with its canonical orientation, i.e. $\int_{\delta\Gamma}\frac{dz}{z}=2\pi i$. We define:
$$\tau_i(A)=(-1)^{i}\delta\eta_i(A)\cdot \Gamma\cdot (\delta\Gamma)^{i-1}+\eta_i(A)\cdot(\delta\Gamma)^i$$
where $\cdot$ denotes the usual alternating product (i.e. we take the usual external product to obtain a cycle
inside $\Box^{2n}$ and then take the alteration with respect to the action of the group $\Sigma_{2n}\rtimes (\Z/2\Z)^{2n}$). And $\delta$
denotes the topological boundary defined for a cycle
$f(s_1,\dots,s_n)$ by:
\begin{align*}
(\delta f)(s_1,\dots
s_{n-1})=f(0,s_1,s_2,\dots,s_{n-1})-f(s_1,s_1,s_2,\dots,s_{n-1})+\cdots
\\
 \cdots+(-1)^{n-1}
f(s_1,\dots,s_{n-1},s_{n-1})+(-1)^{n}f(s_1,\dots,s_{n-1},1).
\end{align*}
Finally denote
$$\xi_{\gamma}(A)=\sum_{k=1}^n \tau_k(A).$$
Note that $\xi_{\gamma}(A)\in \tilde{\D}^0(n)$, i.e. it is a topological cycle inside $\Box_{\C}^{2n}$ of (real) dimension $2n$.  We will denote the corresponding element in the following homological version of $\tilde{\D}$ also by the same notation.
$$\D^0(n) = \mbox{Alt}\varinjlim\h_{2n}(S\cup{\mathcal J}^{2n},{\mathcal J}^{2n})$$
(where ${\mathcal J}^{2n}$ is the union of all the codimension 1 hyper
planes  of $({\Bbb P}^1)^{2n}$ obtained by letting one coordinate equal
to 1. The limit is taken over certain admissible subsets of $\C^{2n}$ of (real) dimension $2n$.)

We define  $\D'$  to be the sub DGA of topological (homological) cycles generated by  $\sigma(\rho(A_i))$ and $\xi_{\gamma}(A_i)$ where $A_1,\dots,A_n$ is any cut of $A$ and $\sigma :\N\To\D$ is the underlying topological cycle (i.e. the fundamental class). In fact one can take the algebra generated by these elements. It will automatically be a DGA because of the following proposition:
\begin{prop}\label{diff}
The differential of $\xi_{\gamma}(A)$ is given by the
formula
\begin{equation*}
\sigma(\rho(A))+\sum \xi_{\gamma}(A')\sigma(\rho(A''))
\end{equation*}
where we sum over all 2-cuts $(A',A'')$ of $A$.
\end{prop}
{\bf Proof.}  Since $\eta_k$ were defined in terms of $\rho_i$'s, they satisfy:
$$d\eta_k(A)=-\delta\eta_{k+1}(A)+(-1)^k\sum \eta_k(A')\sigma(\rho(A''))$$
$$d(\delta\eta_k(A))=(-1)^k\sum \delta\eta_k(A')\sigma(\rho(A''))$$
 Hence:
\begin{align*}
d\tau_k(A)=& (-1)^{k}d(\delta \eta_k(A)\cdot \Gamma\cdot
\delta\Gamma^{k-1})+ d\eta_k(A)\cdot \delta\Gamma^k
\\
=& (-1)^kd(\delta\eta_k(A))\cdot \Gamma\cdot \delta\Gamma^{k-1}-\delta\eta_k(A)d(\Gamma)\delta\Gamma^{k-1}+
d\eta_k(A)\cdot \delta\Gamma^k
\\
=&\sum_{\small{\mbox{2-cuts of A}}}\Bigl\{\delta\eta_k(A')\sigma(\rho(A''))\cdot\Gamma\cdot\delta\Gamma^{k-1}
+(-1)^k\eta_k(A')\sigma(\rho(A''))\cdot\delta\Gamma^k\Bigr\}
\\
&\:\:\:\:+\delta\eta_k(A)\cdot\delta\Gamma^{k-1}-\delta\eta_{k+1}(A)\cdot\delta\Gamma^k
\\
=&\sum_{\small{\mbox{2-cuts of A}}}\Bigl\{(-1)^k\delta\eta_k(A')\cdot\Gamma\cdot\delta\Gamma^{k-1}
+\eta_k(A')\cdot\delta\Gamma^k\Bigr\}\cdot\sigma(\rho(A'')) \\
&\:\:\:\:+\delta\eta_k(A)\cdot\delta\Gamma^{k-1}-\delta\eta_{k+1}(A)\cdot\delta\Gamma^k
\\
=&\delta\eta_k(A)\cdot\delta\Gamma^{k-1}-\delta\eta_{k+1}(A)\cdot\delta\Gamma^k
+ \sum_{\small{\mbox{2-cuts of A}}} \tau_k(A')\sigma(\rho(A'')) \end{align*}
Summing over $k$ we get:
$$d\xi_{\gamma}(A)=\delta\eta_1(A)+\sum_{\small{\mbox{2-cuts of A}}}\xi_{\gamma}
(A')\cdot\sigma(\rho(A''))$$
but since $\delta\eta_1(A)=\sigma(\rho(A))$ the proposition is proved.\qed

Finally our last piece of data, is a construction of certain element $Z_{\gamma}(A)\in \h^0B(\D',\N')$.

\begin{defn}\label{Z} Define $Z_{\gamma}(A)\in \h^0B(\D',\N')^0$ by
$$\sum 1\otimes [\rho(A_1)|\dots|\rho(A_k)]+\sum\xi_{\gamma}(A_1)[\rho(A_2)|\dots|\rho(A_k)]$$
where the sums are over all cuts of $A$.
\end{defn}
\begin{lem}\label{dZ} $dZ_{\gamma}(A)=0$, hence it defines an element of $H^0B(\D',\N')$ denoted by the same notation.
\end{lem}
{\bf Proof.} The first sum in $Z_{\gamma}$ is in fact $1\otimes \I(A)$. Since $d\I(A)=0$ as an element of $B(\N')^0$, the only contribution of the differential for this term is
$$-\sum \sigma(\rho(A_1))\otimes [\rho(A_2)|\cdots|\rho(A_k)].$$
The external differential of each summand of the second sum is
$$d(\xi_{\gamma}(A_1))\otimes[\rho(A_2)|\cdots|\rho(A_k)]-\sum_{i=2}^n \xi_{\gamma}(A)\otimes[\rho(A_1)|\cdots|d\rho(A_i)|\cdots|\rho(A_k)]$$
Since $d\xi_{\gamma}(A_1)=\sigma(\rho(A_1))+\sum \xi_{\gamma}(A_1')\sigma(\rho(A_1''))$ and $d\rho(A_i)=-\sum \rho(A_i')\rho(A_i'')$, adding the above two terms gives us:
\begin{align*}
\sum_{\small{\mbox {2-cuts of $A_1$}}}\xi_{\gamma}(A_1')\sigma(\rho(A_1''))&\otimes [\rho(A_2)|\cdots|\rho(A_k)]
\\ +\sum_{i=2}^n \xi_{\gamma}(A)\otimes[\rho(A_1)|&\cdots|\sum\rho(A_i')\rho(A''_i)|\cdots|\rho(A_k)]
\end{align*}
If we only sum this over $k$-cuts of $A$, it will be canceled by the internal differential of
$$\sum_{\mbox{\small{$(k+1)$-cuts}}}\xi_{\gamma}(A_1)\otimes [\rho(A_2)|\cdots|\rho(A_{k+1})].$$\qed

\begin{lem}\label{lambda}
Assume $a_0\ne a_1$ and $a_n\ne a_{n+1}$. Under the map (see lemma \ref{BK})
$$\Lambda:H^0B(\D',\N')\To H^0B(\N')\otimes \C$$
$Z_{\gamma}(A)$ maps to
$$\I(A)+ \sum \left(\int_{\gamma}\omega_{A_1}\right)[\rho(A_2)|\cdots|\rho(A_k)]$$ where the sum is taken over all cuts $A_1,\dots,A_k$ of $A$, and if $A_1=(a_0; a_{i_1},\dots, a_{i_m}; a_{n+1})$
$$\omega_{A_1}=(2\pi i)^{-m}\frac{dt}{t-a_{i_1}}\circ\dots\circ\frac{dt}{t-a_{i_m}}.$$
\end{lem}
{\bf Proof.} It is enough to show that
$$\lambda(\xi_{\gamma}(A_1))=x^m\int_{\gamma}\omega_{A_1}.$$
where $m=\deg(A_1)$. To show this, first observe that if $i<m$
$$\int_{\eta_i(A_1)}\frac{dz_1}{z_1}\wedge \cdots\wedge \frac{dz_{2m-i}}{z_{2m-i}} =0$$
This is so, since at least two components have the same parameter, and the alternation will kill them. By a similar argument we can, replace $\eta_i(A)$ with $\delta\eta_i$ and this time we can even include $i=m$.(This is the only place that one needs to look inside the machinery of integration theory). So the only term
in $\xi_{\gamma}(A_1)$ that gives a nonzero contribution is $\eta_m(A_1)\cdot (\delta\Gamma)^m$. Since $\int_{\delta\Gamma}\frac{dz}{z}=2\pi i$ we get:
$$\lambda(\xi_{\gamma}(A_1))=x^m (2\pi i)^{-m}\int_{\eta_m(A_1)}\frac{dz_1}{z_1}\wedge\dots\wedge\frac{dz_m}{z_m}.$$
It is left to check that:
$$\int_{\eta_m(A_1)}\frac{dz_1}{z_1}\wedge\dots\wedge\frac{dz_m}{z_m}=\int_{\gamma}\frac{dt}{t-a_{i_1}}\circ\dots\circ\frac{dt}{t-a_{i_m}}$$
But $\eta_m(A_1)$ is the cycle
$$\mbox{Alt}(c_1(\gamma(s_1)-a_{i_1}),\dots,c_m(\gamma(s_m)-a_{i_m}))$$ where $0\le s_1\dots\le s_m\le 1$ and the constants $c_i$ depend on our integration theory. However since the form
$\frac{dz_1}{z_1}\wedge \dots \wedge \frac{dz_m}{z_m}$ is invariant under the action of $(\C^*)^m$ and is alternating under the action of $(\Z/2)^m\rtimes \Sigma_m$, we get the iterated integral on the right, by definition.

Finally let us show that the formula for $\Lambda(Z_{\gamma}(A))$ which a priory is in $B(\N)\otimes \C$, is in fact inside $\h^0B(\N)\otimes \C$. This is not needed because $\Lambda$ is map of differential graded vector spaces, however
it is not hard to see it directly. In fact if one fixes $A_1=(a_0;a_{i_1},\dots,a_{i_m};a_{n+1})$ then $A_2,\dots,A_k$ are a shuffle of cuts for $(a_0;\dots;a_{i_1})$, $(a_{i_1};\dots;a_{i_2}), \dots, (a_{i_m};\dots;a_{n+1})$ and so the formula can be written as:
$$\I(A)+ \sum \left(\int_{\gamma}\omega_{A_1}\right)\I(A_2)\cdots\I(A_k)$$
where this time we only sum over elementary cuts of $A$.
 \qed

\begin{lem} If $a_0\ne a_1$ and $a_n\ne a_{n+1}$, the pair $(\N',\D')$ satisfies the conditions in definition \ref{admissible} of appendix B.
\end{lem}
{\bf Proof.} Condition (1) follows from lemma \ref{generator}. Condition (2) is trivial by definitions. Condition (3) follows from the calculations of the previous lemma. Conditions (4) and (6) are left to the reader to verify. To check condition (5) we need to show that $\tau^*(\I(A_i))=0$ where $A_i$ is a cut component of $A$. However:
$$\Z_{\gamma}(A_i)= 1\otimes \I(A_i)+\sum \xi_{\gamma}(A_{i1})[\rho(A_{i2})|\dots|\rho(A_{ik})]$$
(summing over all cuts of $A_i$) has differential zero and $d(1\otimes \I(A_i))=-\tau(\I(A_i))$ so
$$\tau(\I(A_i))= d(\sum \xi_{\gamma}(A_{i1})[\rho(A_{i2})|\dots|\rho(A_{ik})]).$$
\qed

\begin{thm}\label{compare}
Let $A=(a_0;a_1,\dots,a_n;a_{n+1})$.
The Hodge realization of ${\mathbb I}(A)$ is $\I^{\mathcal H}(A)$.
\end{thm}
{\bf Proof.} According to the lemma \ref{convergent} for a permuting integration theory we can assume, $a_0\ne a_1$ and $a_n\ne a_{n+1}$. Recall the MHTS's $\Pi^n({\mathbb A}^1-\{a_1,\dots,a_n\};a_0,a_{n+1})$ and $H(\N',n)[n]$ of appendix B. Their $(-2k)$-de Rham graded pieces are homogeneous degree $k$ elements of $\Q\langle\langle X_{a_i}\rangle\rangle$ and $\h^0B(\N')(n-k)$ respectively. We define the following map from the first space to the second:
$$X_{b_1}\cdots X_{b_k}\mapsto \sum[\rho(A_2)|\cdots|\rho(A_m)]$$
where we sum over all cuts $A_1,\dots,A_m$ of $A$ such that  $A_1=(a_0;b_1,\dots,b_k;a_{n+1})$. If there are no such cuts the value is zero and if $m=1$, the value is 1. Using the same ideas as in the lemma \ref{lambda}, one can see that this in fact lands in $\h^0B(\N')(n-k)$. It maps $1$ to $\I(A)$. It also maps $X_{a_1}\dots X_{a_n}$ to $1$, while mapping any other monomial of degree $n$ to zero. Hence it does respect the frames of $\I^{\mathcal H}(A)$ and $\mbox{Real}_{MHTS}(\I(A))$. It remains to show that it is a morphism of MHTS's. i.e. it maps the Betti subspace of $\C\langle\langle X_{a_i}\rangle\rangle$ defined by the image of
$$\Phi:\Pi^N(X,a_0,a_{n+1})\To\C\langle\langle X_{a_i}\rangle\rangle$$
$$\gamma\mapsto  1+\sum_{b_i\in\{a_1,\dots,a_n\}} (2\pi i)^{-k}(\int_{\gamma}\frac{dt}{t-b_1}\circ\dots\circ\frac{dt}{t-b_k}) X_{b_1}\dots X_{b_k}$$
to the Betti subspace of $\h^0B(\N')\otimes \C$ defined by the image of the map $\Lambda$ of lemma \ref{BK}.
 Observe that under the map defined between the two structures above, $\Phi(\gamma)$ maps to
$$\I(A) + \sum (\int_{\gamma}\omega_{A_1})[\rho(A_2)|\dots|\rho(A_m)].$$
where we sum over all cuts of $A$. This is $\Lambda(Z_{\gamma}(A))$ and therefore is inside the Betti subspace of $H(\N',n)$.\qed
\appendix
\section{Bar Construction}

In this section we recall the basic definitions of bar construction. We follow the exposition of Hain and Zucker \cite{HZ}. It differs from \cite{BK} in certain signs.

A differential graded algebra (DGA for short) over a field $F$ is a graded vector space over $F$:
$${\mathcal A}=\bigoplus_{n\in \Z}{\mathcal A}^n$$ with an associative multiplication and linear maps:
$$d:{\mathcal A}^n\To {\mathcal A}^{n+1}$$
such that $d\circ d = 0$. For $a\in{\mathcal A}^n$ and $b\in{\mathcal A}^m$, $a\cdot b\in {\mathcal A}^{n+m}$ and:
$$d(a\cdot b)=(da)\cdot b+(-1)^n a\cdot(db).$$
An augmentation for $\mathcal A$ is an onto map $\epsilon:{\mathcal A}\To F$ of DGA's, where $F$ is the trivial DGA, concentrated in degree $0$. If we assume that $\mathcal A$ has a unit, we also require that $\epsilon(1)=1$.

A (right) module over $\mathcal A$ is a graded vector space $\mathcal M=\oplus_n {\mathcal M}^n$ with a differential $d:{\mathcal M}^n\To {\mathcal M}^{n+1}$ and an action:
$${\mathcal M}\otimes {\mathcal A}\To {\mathcal M}$$
such that:
\begin{enumerate}
\item For $m\in\mathcal M$ and $a,b\in \mathcal A$, $(ma)b=m(ab)$.
\item If $\mathcal A$ has a unit the $m\cdot 1 = m$.
\item For $m\in {\mathcal M}^n$ and $a\in {\mathcal A}^m$, $am\in {\mathcal M}^{n+m}$ and
$$d(ma)=(dm)a+(-1)^nm(da).$$
\end{enumerate}

One can similarly define a left $\mathcal A$ module. We now give the details of the bar construction. Let $\mathcal A$ be an augmented DGA and ${\mathcal A}^+$ be the kernel of the augmentation. Let $\mathcal M$ be a right $\mathcal A$ module. The bar construction assigns a differential graded vector space $B({\mathcal M},{\mathcal A})$ to this data as follows. As a graded vector space it is:
$${\mathcal M}\otimes_F T({\mathcal A}^+[1]).$$
Here $T({\mathcal A}^+[1])$ is the tensor algebra over the shifted graded vector space ${\mathcal A}^+[1]$. Hence the degree of a homogeneous element $m\otimes a_1\otimes \dots \otimes a_r$ is
$$\deg (m)+\deg(a_1)+\dots+\deg(a_r)-r.$$
Often we use the bar notation: $m\otimes[a_1|\dots|a_r]$ for such an element. There are two differential on this space, the external differential denoted by $d_{ext}$ and the internal differential denoted by $d_{int}$. They are defined by the following formulae:
\begin{align*}
d_{ext}(m\otimes[a_1|\dots |a_r])=&dm\otimes[a_1|\dots|a_r] \\
&+\sum_{i=1}^r Jm\otimes[Ja_1|\dots|Ja_{i-1}|da_i|a_{i+1}|\dots|a_r],
\\
d_{int}(m\otimes[a_1|\dots|a_r])=&(Jm\cdot a_1)\otimes[a_2|\dots|a_r]
\\
&+\sum_{i=1}^{r-1} Jm\otimes [Ja_1|\dots|Ja_{i-1}|Ja_{i}\cdot
a_{i+1}|\dots|a_r].
\end{align*}
 the operation $J$ is given by $J(a)=(-1)^{\deg(a)-1}a$ on
the homogeneous elements. Note that $\deg(a)-1$ is the degree of $a$
in the shifted complex $\mathcal A[1]$. If $\mathcal M = F$ is the trivial DGA, instead of $B(F,\mathcal A)$ we write $B(\mathcal A)$ and we write $[a_1|\dots|a_k]$ for $1\otimes [a_1|\dots|a_k]$.

There is a coproduct on $B(\mathcal A)$ (where the empty tensor is 1 by convention):
$$\Delta([a_1|\dots|a_r])=\sum_{s=0}^r [a_1|\dots|a_s]\otimes
[a_{s+1}|\dots|a_r].$$

Up to this point all the constructions work for a DGA which is not
necessarily commutative. The product is defined only when both $\mathcal M$
and $\mathcal A$ are commutative DGA's (i.e. $ab=(-1)^{nm}ba$ for $a\in \mathcal A^n$ and $b\in \mathcal A^m$) by the shuffle:
$$m\otimes[a_1|\dots|a_r]\cdot m'\otimes[a_{r+1}|\dots|a_{r+s}]:=\sum
sgn_{{\bf{a}}}(\sigma)(m\cdot
m')\otimes[a_{\sigma(1)}|\dots|a_{\sigma(r+s)}]$$ where $\sigma$ runs over
the $(r,s)$-shuffles and the sign is obtained
by giving $a_i$'s weights $\deg(a_i)-1$. For example if all $a_i$'s are of odd weights then all the signs are plus.
These definitions make $B({\mathcal A})$ into a Hopf algebra where $\mathcal A$ is a commutative DGA.

\section{Mixed Hodge-Tate Structures}
We recall a definition for the mixed Hodge-Tate structures (MHTS for short) below:
\begin{defn} A MHTS $H$ is a finite dimensional $2\Z$-graded $\Q$-vector
space $H_{dR}=\oplus_n H_{2n}$ together with a $\Q$-subspace $H_B$
of $H_{dR}\otimes \C$ such that\footnote{It is more conventional to replace the right hand side with $(2\pi i)^{-m}H_{2m}$, however our definition is equivalent and make our constructions a little easier.} for all $m$:
$$\mbox{Im}\left(H_B\cap (\bigoplus_{n\le m} H_{2n}\otimes \C)\stackrel{Proj}\To
H_{2m}\otimes \C\right)= H_{2m}.$$ A morphism $f$ is a graded
morphism $f_{dR}$ from $H_{dR}$ to $H_{dR}'$ such that
$f_{dR}\otimes 1:H_{dR}\otimes \C\To H_{dR}'\otimes \C$ sends $H_B$
to $H_B'$, this induced map is denoted by $f_B$.
\end{defn}
\begin{rem}\label{shift} If $H=(H_{dR},H_B)$ is a MHTS as above, the shifted structure $H[n]$ is given by
$$H[n]_{dR}=\bigoplus_k H[n]_{2k}$$
$$H[n]_{2k}=H_{2k+2n}$$
and the same $H_B$.
\end{rem}
Our two basic examples of such structures are the truncated torsor of paths on $X={\Bbb A}^1\backslash\{a_1,\dots,a_n\}$
$$\Pi^N(X;a_0,a_{n+1})$$
and the space $\h^0B(\N')$ for certain sub-DGA's $\N'$ of $\N$. Let us recall the construction of MHTS on these two spaces. For the first one we follow \cite{G} \S 5 and for the second one we follow \cite{BK} \S 8.

Assume that $a_1,\dots, a_n$ are distinct. Let $H_{-2m}$ be the space of homogeneous polynomials of degree $m$ on non-commuting variables $X_1,\dots,X_n$ with coefficients in $\Q$. This is a space of dimension $m^n$ for $m>0$ and $1$ for $m=0$ over $\Q$. Let
$$H_{dR}=\bigoplus_{0\le m \le N}H_{-2m}.$$
 Let $\omega_i:=(2\pi i)^{-1}d\log(t-a_i)$. Any path $\gamma:[0,1]\To X$ from $a_0$ to $a_{n+1}$ defines an element of $H_{dR}\otimes \C$ by the Feynman-Dyson formula:
$$\Phi(\gamma)=\sum \left(\int_{\gamma}\omega_{i_1}\circ\dots\circ\omega_{i_k}\right)X_{i_1}\dots X_{i_k}$$
where we sum over all $(i_1,\dots,i_k)$ for $0\le k\le N$ and $i_j\in\{1,\dots,n\}$. One can show that this element depends only on homotopy class of $\gamma$ as a path from $a_0$ to $a_{n+1}$ in $X$, and if we extend $\Phi$ by linearity to $\Q[\pi(X;a_0,a_{n+1})]$ it vanishes on $I_{a_0}^{N+1}\Q[\pi(X;a_0,a_{n+1})]$. $I_{a_0}$ is the augmentation ideal $\mbox{Ker}(\Q[\pi_1(X;a_0)]\To \Q)$. Hence we get a map
$$\Phi: \Pi^N(X;a_0,a_{n+1})\To H_{dR}\otimes \C.$$
Its image is the Betti subspace $H_B$.

If $a_0$ or $a_{n+1}$ are in $\{a_1,\dots, a_n\}$ one needs to consider tangential base point and regularization of the iterated integrals as follows. Identify the tangent space at each point of $\C$ with $\C$ and as the tangential base point at $a$, take the unit tangent vector $1$. Now by a path from $a_0$ to $a_{n+1}$ we mean a differentiable map $\gamma :[0,1]\To \C$ such that $\gamma(0)=a_0$, $\gamma(1)=a_{n+1}$, $\gamma'(0)=\gamma'(1)=1$ and $\gamma$ restricted to the open interval $(0,1)$ maps to $\C-\{a_1,\dots,a_n\}$. The regularized iterated integral is defined as follows. Let $\gamma^*\omega_i=f_i(t)dt$, then the integral
$$\int_{\epsilon\le t_1\le \dots\le t_k\le 1-\epsilon}f_{i_1}(t_1)\dots f_{i_k}(t_k)dt_1\dots dt_k$$
can be written as $\sum_{i=1}^m c_m(\epsilon)\log^i \epsilon$ where $c_m$ are holomorphic around 0. The regularized value is $c_0(0)$.
\begin{lem}
For a given $N$, the $2\Z$-graded vector space
$$H_{dR}=\oplus_{k\le N} H_{-2k}$$
where $H_{-2k}$  is the space of homogeneous polynomials of degree $k$ on non-commuting variables $X_1,\dots,X_n$ with coefficients in $\Q$, together with the $\Q$ subspace $H_B$ of $H_{dR}\otimes \C$ defined by:
$$\mbox{Im}(\Phi:\Pi^N(X;a_0,a_{n+1})\To H_{dR}\otimes \C)$$
is a MHTS. It is also denoted by $\Pi^N(X;a_0,a_{n+1})$.
\end{lem}

Now we recall the MHTS defined in \cite{BK} on $\h^0B(\N')$ for certain sub-DGA's $\N'$ of $\N$.

Let an admissible topological cycle of (real) codimension $2r$ in $\Box_{\C}^n$ be a linear sum of admissible subsets
 of $\Box_{\C}^n$ which are a union of smooth singular simplices of (real) codimension $2r$ and disjoint interiors. Admissible means that the intersection with any face of $\Box^n$ is either empty is also a union of singular simplices of real codimension $2r$ and disjoint interiors. This group is denoted by $\mbox{TCycle}^r(n)$. Let
$${\tilde{\D}}^n(r)=\mbox{Alt}(\mbox{TCycle}^r(2r-n)\otimes \Q)$$
A similar constructions to $\N$, will make $\tilde{\D}=\oplus \tilde{\D}^n(r)$ into a DGA with a natural morphism:
$$\sigma:\N\To \tilde{\D}$$
that assigns the underlying topological cycle to an algebraic cycle. We would like to consider these topological cycles up to homology. i.e. If for $Z\in \tilde{\D}^0(n)$ we have $dZ=0$, we want it to be equivalent to $(\delta\Gamma)^{2n}$, where $\Gamma$ is a small disk around zero in $\Box_\C$ and $\delta$ is the boundary. This is done in more details in \cite{BK}. But the idea is to consider
$${\D}^n(r):=\mbox{Alt}\varinjlim\h_{2r-2n}(S\cup{\mathcal
J}^{2r-n},{\mathcal J}^{2r-n})$$
where ${\mathcal J}^n$ is the union of all the codimension 1 hyper
planes  of $({\Bbb P}^1)^n$ obtained by letting one coordinate equal
to 1. The limit is taken over all admissible subset $S$ in ${\C}^{2r-n}$ of (real) codimension $2r$. This has a
natural structure of DGA (refer to \cite{BK}, \S 8).
\\
\begin{defn}\label{admissible}
A pair of a sub-DGA $\N'$ of $\N$ and a sub-DGA $\D'$ of $\D$ is called admissible if:
\begin{enumerate}
\item Each space $\h^0B(\N')(r)$ is finite dimensional $\Q$-vector space.
\item The map $\sigma$ maps $\N'$ to $\D'$.
\item Let $\C[x]$ be a DGA concentrated in degree zero with Adams grading given by powers of $x$ and zero differential. Define $\lambda:\D'\To \C[x]$  to be zero on $(\D')^n(r)$ for $n\ne 0$ and for $c\in (\D')^0(r)$ , i.e. an admissible topological cycle in $\Box^{2r}$ of real dimension $r$:
    $$\lambda(c)=\left((2\pi i)^{-2r}\int_c \frac{dz_1}{z_1}\wedge\dots\wedge\frac{dz_{2r}}{z_{2r}}\right) x^r$$
    We will assume that $\lambda$ is a well-defined morphism of DGA's. We call $\lambda$ the period map.
\item The map $\lambda:\h^0(\D')\To \h^0(\C[x])=\C[x]$ defines an isomorphism $\h^0(\D')=\Q[x]$. (This condition is missing from \cite{BK}, however is needed if one wants to prove that the following structure is a MHTS)
\item Let $\tau: B(\N')\to B(\D',\N')[1]$ (bracket denotes the shift in the degrees) be given by:
$$[a_1|\dots|a_r]\mapsto \sigma(a_1)\otimes[a_2|\dots|a_n].$$
It can easily be seen that this is a map of differential vector spaces so it induces:
$$\tau^* : \h^0B(\N')\To\h^1B(\D',\N').$$
We will assume that this map is zero.
\item There is one more condition that is technical and we refer the reader to page 597 of \cite{BK} for its statement.
\end{enumerate}
\end{defn}
Define
$$\Lambda: \h^0B(\D',\N')\To\h^0B(\N')\otimes \C$$
via the period map followed by evaluation at $x=1$:  $\D'\stackrel{\lambda}\To\C[x]\stackrel{x=1}\To \C$.
\begin{lem}\label{BK}
For a given $n$, the $2\Z$-graded vector space:
$$H_{dR}=\oplus_{r\le n} H_{2r},\:\: H_{2r}=H^0B(\N')(r)$$
together with the $\Q$-subspace $H_B$ of $H_{dR}\otimes \C$ given by:
$$\mbox{Im}(\Lambda: \oplus_{r\le n}H^0B(\D',\N')(r)\To \oplus_{r\le n}H^0B(\N')(r)\otimes \C)$$
is a MHTS. We denote this MHTS by $H(\N',n)$
\end{lem}

The category of MHTS's is a neutral Tannakian category with a fiber functor given by sending $H$ to $H_{dR}$. The general formalism of Tannakian categories imply that this category is equivalent to the category of graded comodules over a certain Hopf algebra. Beilinson observed that this Hopf algebra can be described using framed MHTS. An $n$-framed MHTS is a MHTS, $H$ together with elements $v\in H_0$ and $\hat{v}\in \hm(H_{-2n},\Q)$. Two $n$-framed MHTS's, $(H,v,\hat{v})$ and $(H',v',\hat{v}')$ are equivalent if there is a morphism $f:H\To H'$ such that $f_{dR,0}(v)=v'$ and $\hat{v}'\circ f_{dR,-2n} =\hat{v}$. This generates an equivalence relation, whose classes is denoted by ${\chi}_{MHTS}(n)$.
The space ${\chi}_{MHTS}=\oplus_{n\ge 0} {\chi}_{MHTS}(n)$ has a natural structure of a Hopf algebra. The realization functor
$$Real_{MHTS}:\h^0B(\N)\To \chi_{MHTS}$$
as a morphism of Hopf algebras has been defined in section 7 of \cite{BK}. It is very complicated, however a more
accessible definition at least for a sub-DGA, $\N'$ and a relative sub-DGA $\D'$ satisfying the above assumptions is also given.
\begin{defn}
The Hodge realization map $\mbox{Real}_{MHTS}:\h^0B(\N')\To\chi_{MHTS}$ is given by sending $a\in \h^0B(\N')(n)$
to the shifted MHTS $H(\N',n)[n]$ (ref. to the remark \ref{shift}) with the following frames:
$$a\in \h^0B(\N')(n)=H(\N',n)[n]_0$$
$$\epsilon\in\hm(\h^0B(\N')(0),\Q)=\hm(H(\N',n)[n]_{-2n},\Q)$$
where $\epsilon$ is the augmentation.
\end{defn}

Another important class of framed MHTS are associated to the iterated integrals.
\begin{defn}\label{framed integral}
For a sequence $(a_0;a_1,\dots,a_n;a_{n+1})$ of elements of $\C$, the framed MHTS $\I^{\mathcal H}(a_0;a_1,\dots,a_n;a_{n+1})$ is defined by the MHTS $H=\Pi^n(\C-\{a_1,\dots,a_n\};a_0,a_{n+1})$ with frames given by
$$1\in H_0=\Q$$
$$(X_{a_1}\cdots X_{a_n})'\in \hm(H_{-2n},\Q)=\hm(\Q\langle\langle X_{a_i}\rangle\rangle_{\deg =n},\Q)$$
where $(X_{a_1}\cdots X_{a_n})'$ is a morphism that sends $X_{a_1}\cdots X_{a_n}$ to $1$ and any other monomial
$X_{b_1}\cdots X_{b_n}$ with $b_i\in\{a_1,\dots,a_n\}$ to zero.
\end{defn}

\end{document}